\documentclass[12pt]{article}

\usepackage{amsmath,amsthm}
\usepackage{amssymb}
\usepackage{color}
\usepackage{yfonts}
\usepackage{mathrsfs,textcomp}
\usepackage{enumerate}

\usepackage{latexsym}
\usepackage[english]{babel}
\usepackage{graphicx}
\usepackage{verbatim}
\renewcommand{\exp}[1]{\mathrm{exp}\Pt{#1}}
\newcommand{\der}[2][]{\frac{\partial #1}{\partial #2}}
\newcommand{\G}{G}
\newcommand{\GG}{\mathfrak{G}}
\renewcommand{\L}{\mathfrak{L}}
\renewcommand{\l}{\mathfrak{l}}
\newcommand{\p}{\mathfrak{p}}
\renewcommand{\P}{\mathfrak{P}}
\newcommand{\gruppi}{$\GG_4$ and $\GG_5$}
\newcommand{\Hp}{{\bf (H$ _0$)}}
\newcommand{\Repbase}{\mathfrak{X}}
\newcommand{\dRepbase}{d\Repbase}
\newcommand{\domRbase}{{\cal H}}
\newcommand{\fzRbase}{\psi}
\newcommand{\eigRbase}{\al}
\newcommand{\DHLbase}{\hat{\Delta}_\SR}

\renewcommand{\b}[1]{{\bf #1}}
\newcommand{\SR}{H}

\newcommand{\Dh}{\Delta_H}
\renewcommand{\span}[1]{\mathrm{span}\Pg{#1}}
\newcommand{\cont}{{\cal C}}
\newcommand{\distr}{{\blacktriangle}}

\newcommand{\hp}{hypothesis}

\renewcommand{\H}{{\cal H}}
\newcommand{\g}{{\mathbf g}}
\newcommand{\metr}{\g}
\renewcommand{\div}{\mathrm{div}}

\newcommand{\gradh}{\mathrm{grad}_\SR}

\newcommand{\Mat}[2]{\Pt{\ba{#1} #2\ea}}

\newcommand{\funz}[5]{#1 : \begin{tabular}{ccl}
 #2 &$\rightarrow$& #3 \\
 #4 & $\mapsto$& #5   \end{tabular}}

\newcommand{\mmfunz}[5]{\begin{center}
\funz{$\displaystyle{#1}$}{$\displaystyle{#2}$}{$\displaystyle{#3}$}{$\displaystyle{#4}$}{$\displaystyle{#5}$}
\end{center}}
\newcommand{\fz}[3]{#1:\, #2 \rightarrow #3}

\newcommand{\bbibitem}{\bibitem}

\newcommand{\llabel}[1]{{\label{#1}}}

\newcommand{\ffoot}[1]{}
\newcommand{\fffoot}[1]{}

\renewcommand{\r}[1]{(\ref{#1})}

\newcommand{\bi}{\begin{itemize}}
\newcommand{\ei}{\end{itemize}}
\newcommand{\be}[1][]{\begin{enumerate}[#1]}
\newcommand{\ee}{\end{enumerate}}

\newcommand{\bd}{\begin{description}}
\newcommand{\ed}{\end{description}}
\renewcommand{\i}{\item}

\newcommand{\bqn}{\begin{eqnarray}}
\newcommand{\eqn}{\end{eqnarray}}
\newcommand{\eqnn}{\nonumber\end{eqnarray}}
\newcommand{\eqnl}[1]{\llabel{#1}\end{eqnarray}}

\newcommand{\nn}{\nonumber\\}
\newcommand{\ba}[1]{\begin{array}{#1}}
\newcommand{\ea}{\end{array}}

\newcommand{\R}{\mathbb{R}}
\newcommand{\C}{\mathbb{C}}
\newcommand{\N}{\mathbb{N}}

\newcommand{\bproof}{\begin{proof}}
\newcommand{\eproof}{\end{proof}}
\newtheorem{Theorem}{\bf Theorem}
\newtheorem{lemma}[Theorem]{\bf Lemma}
\newtheorem{corollary}[Theorem]{\bf Corollary}
\newtheorem{definition}[Theorem]{\bf Definition}
\newtheorem{proposition}[Theorem]{\bf Proposition}
\newtheorem{remark}[Theorem]{\bf %\underline
{{\sl Remark}}}

\newcommand{\bt}{\begin{Theorem}}
\newcommand{\et}{\end{Theorem}}
\newcommand{\bl}{\begin{lemma}}
\newcommand{\el}{\end{lemma}}
\newcommand{\bp}{\begin{proposition}}
\newcommand{\ep}{\end{proposition}}
\newcommand{\bc}{\begin{corollary}}
\newcommand{\ec}{\end{corollary}}
\newcommand{\bdeff}{\begin{definition}}
\newcommand{\edeff}{\end{definition}}
\newcommand{\brem}{\begin{remark}\rm}
\newcommand{\erem}{\end{remark}}
\newcommand{\auth}[1]{{\sc #1}}
\newcommand{\tit}[1]{{\rm #1}}
\newcommand{\titl}[1]{{\it #1}}
\newcommand{\jou}[1]{{\it #1}}
\newcommand{\vol}[1]{{\it #1}}
\newcommand{\pp}[1]{pp.~#1}
\newcommand{\lam}{\lambda}
\newcommand{\al}{\alpha}

\newcommand{\ga}{\gamma}
\newcommand{\de}{\delta}

\renewcommand{\th}{\theta}
\newcommand{\Id}{\mathrm{Id}}
\newcommand{\Lie}{\mathrm{Lie}}
\newcommand{\Pt}[1]{\left( #1 \right)}
\newcommand{\Pg}[1]{\left\{ #1 \right\}}
\newcommand{\Pq}[1]{\left[ #1 \right] }
\newcommand{\Pa}[1]{\langle #1 \rangle}
\newcommand{\GFT}{GFT}

\newcommand{\FF}{{\cal F}}
\renewcommand{\div}{\mbox{div}}

\newcommand{\HS}{{\mathbf{HS}}}

\newcommand{\gr}[1]{{#1}}
\begin{document}

\thispagestyle{empty}
\begin{center} \noindent
{\LARGE{\sl{\bf Hypoelliptic heat kernel over 3-step nilpotent Lie groups}}}
%\vskip 1cm
%\today
\vskip 1cm
Ugo Boscain\footnote{The first author has been supported by the European Research Council (ERC StG 2009 GeCoMethods).}

{\footnotesize CMAP, \'Ecole Polytechnique CNRS, Route de Saclay, 91128 Palaiseau Cedex, France - {\tt ugo.boscain@polytechnique.edu}}\\
\vspace*{.5cm}

Jean-Paul Gauthier\\
{\footnotesize Laboratoire LSIS, Universit\'e de Toulon, France - {\tt gauthier@univ-tln.fr}}

\vspace*{.5cm}

Francesco Rossi\\
{\footnotesize BCAM - Basque Center for Applied Mathematics,

Bizkaia Technology Park, Basque Country, Spain
- {\tt rossi@bcamath.org}}
\end{center}

\begin{quotation}
In this paper we provide explicitly the connection between the hypoelliptic heat kernel for some 3-step sub-Riemannian manifolds and the quartic oscillator. We study the left-invariant sub-Riemannian structure on two nilpotent Lie groups, namely the (2,3,4) group (called the Engel group) and the (2,3,5) group (called the Cartan group or the generalized Dido problem). Our main technique is  noncommutative Fourier analysis that permits to transform the hypoelliptic heat equation into a one dimensional heat equation with a quartic potential.
\end{quotation}

\bigskip

\section{Introduction}

The study of the properties of the heat kernel in a sub-Riemannian manifold  drew an increasing attention since the pioneer work of H\"ormander \cite{hormander}. Since then, many estimates and properties of the kernel in terms of the sub-Riemannian distance have been provided (see \cite{lanconelli-book,boscain-polidoro,folland-stein,rot, varopoulos} and references therein). For some particular structures, it is moreover possible to find explicit expressions of the hypoelliptic heat kernels. In general, this computation can be performed only when the sub-Riemannian structure and the corresponding hypoelliptic heat operator present symmetry properties. For this reason, the most natural choice in this field is to consider invariant operators defined on Lie groups. Results of this kind have been first provided in \cite{gaveau,hul} in the case of the 3D Heisenberg group. Afterwards, other explicit expressions have been found first for 2-step nilpotent free Lie groups (again in \cite{gaveau}) and then for general 2-step nilpotent Lie groups (see \cite{beals-HJ,cyg}). We provide in \cite{nostro-kern} the expressions of heat kernels for 2-step groups that are not nilpotent, namely $SU(2),SO(3),SL(2)$ and the group of rototranslations of the plane $SE(2)$. For other examples, see e.g. \cite{zhu1,zhu2}.

In our paper we present the first results, to our knowledge, about the expression of the hypoelliptic heat kernel on the following 3-step Lie groups. The first one is the Engel group $\GG_4$, that is the nilpotent group with growth vector $(2,3,4)$. Its Lie algebra is $\L_4=\mathrm{span}\Pg{\l_1,\l_2,\l_3,\l_4}$, the  generators of which satisfy
\bqn
&&\Pq{\l_1,\l_2}=\l_3,\ \Pq{\l_1,\l_3}=\l_4,\ \Pq{\l_1,\l_4}=\Pq{\l_2,\l_3}=\Pq{\l_2,\l_4}=\Pq{\l_3,\l_4}=0.
\eqnn 
The second example is the Cartan group $\GG_5$, that is the free nilpotent group  with growth vector $(2,3,5)$. Its Lie algebra is $\L_5=\mathrm{span}\Pg{\l_1,\l_2,\l_3,\l_4,\l_5}$ and  generators satisfy
\bqn
&&\Pq{\l_1,\l_2}=\l_3,\ \Pq{\l_1,\l_3}=\l_4,\ \Pq{\l_2,\l_3}=\l_5,\nn
&&\Pq{\l_1,\l_4}=\Pq{\l_1,\l_5}=\Pq{\l_2,\l_4}=\Pq{\l_2,\l_5}=\Pq{\l_3,\l_4}=
\Pq{\l_3,\l_5}=\Pq{\l_4,\l_5}=0.
\eqnn

In both cases, we consider the heat equation with the so-called intrinsic hypoelliptic Laplacian $\Dh$ (in the sense of \cite{nostro-kern}, see also Section \ref{sez-intr}) of the sub-Riemannian structure for which  $\{g\l_1,g\l_2\}$ ($g$ element of the group) is an orthonormal frame. As it has been proved in \cite{nostro-kern}, since $\GG_4$ and $\GG_5$ are unimodular, then the intrinsic hypoelliptic Laplacian is the sum of the square of the Lie derivative with respect to the vector fields $g\l_1,g\l_2$.
 
One interesting feature of these two sub-Riemannian problems is that they present abnormal minimizers (see \cite{yuri4, yuri5}) and it is known that in both cases $\Delta_H$ is not analytic hypoelliptic \cite{christ}. Hence, for these two examples the  Tr\`eves conjecture\footnote{ We recall that  Tr\`eves conjectured in \cite{treves} that the existence of abnormal minimizers  on a sub-Riemannian manifold is equivalent to the loss of analytic-hypoellipticity of the sub-Laplacian.}  holds. Having information about the expression of the heat kernel can help for further investigations in this direction. 

%These two examples are interesting also for the open questions about analytic-hypoellipticity of the Laplacian. Indeed, Tr\`eves conjectured \cite{treves} that the existence of strict abnormal minimizers of length on a sub-Riemannian manifold is equivalent to the loss of analytic-hypoellipticity of $\Delta_H$. Hence, it is interesting to study the groups \gruppi, since they are examples of sub-Riemannian manifolds in which strict abnormal minimizers appear. The non-analyticity of these heat kernels has been already proved in \cite{christ}.

Any other left-invariant sub-Riemannian structure of rank 2 on these groups is indeed isometric to the ones we study in this paper, see  \cite{yuri4,yuri5}. Moreover, notice that the sub-Riemannian structures we study on $\GG_4$ and $\GG_5$ are local approximations (nilpotentizations, see \cite{gromov}) of arbitrary sub-Riemannian structures at regular points with growth vector $(2,3,4)$ or $(2,3,5)$, hence, roughly speaking, the kernels on $\GG_4$ and $\GG_5$ provide approximations of the heat kernels at these points.\ffoot{un po' naive}\\

The goal of this paper is to transform the hypoelliptic heat equations on these Lie groups into a family of elliptic heat equations on $\R$, depending on one parameter. To this purpose, we apply the method developed in \cite{nostro-kern}, based on the Generalized Fourier Transform (GFT for short), and hence on representation theory of these groups (see \cite[p. 333--338]{dixmier}).

Applying the GFT to the original equation, we get an evolution equation on the Hilbert space where representations act. For both examples, this is the heat equation over $\R$ with quartic potential, the so-called quartic oscillator (see \cite{pham,skala}), for which no general explicit solution is known. Notice that the connection between the quartic oscillator and degenerate elliptic operators has been already noted by previously (see  \cite{greiner-audio}).

It is clearly possible to use numerical approximations of the evolution equation with quartic potential (for which a huge amount of literature is available) to find numerical approximations of the hypoelliptic heat kernel. However, this analysis is outside the aims of this paper.\\

The organization of the paper is the following. In Section \ref{s-Dh} we recall the main definitions from sub-Riemannian geometry, in particular for invariant structures on Lie groups. We then recall the definition of the Generalized Fourier Transform and its main properties. Finally, we recall the main results of our prevuious paper \cite{nostro-kern}, where we studied hypoelliptic heat equations on Lie groups.

The main part of the paper is Section \ref{s-gruppi}. We first present the Lie groups \gruppi, their algebras and their Euclidian and matrix presentations. We then recall  results about their representations. We finally apply the method of computation of hypoelliptic heat kernels to the two groups \gruppi, to  find explicitly the connection between the heat kernels on these groups and the fundamental solution of the 1D heat equation with quartic potential.

\section{The hypoelliptic heat equation on a a sub-Riemannian manifold}
\llabel{s-Dh}
In this section we recall basic definitions from sub-Riemannian geometry, including the one of the intrinsic hypoelliptic Laplacian. Then we recall our method for computing the hypoelliptic heat kernel in the case of unimodular Lie groups, using the GFT.

%we recall an intrinsic definition of the hypoelliptic Laplacian $\Dh$ on a left-invariant sub-Riemannian manifold $M$. We then describe a method for the computation of the hypoelliptic heat kernel on a left-invariant sub-Riemannian manifold. The method is based on the Generalized Fourier Transform, for which we provide definition and properties. For details, refer to our article \cite{nostro-kern}.

\subsection{Sub-Riemannian manifolds}
We start by recalling the definition of sub-Riemannian manifold.
\bdeff
A $(n,m)$-sub-Riemannian manifold is a triple $(M,\distr,{\mathbf g})$, 
where
\bi
\i $M$ is a connected smooth manifold of dimension $n$;
\i $\distr$ is a smooth distribution of constant rank $m< n$ satisfying the {\bf H\"ormander condition}, i.e. $\distr$ is a smooth map that associates to $q\in M$  a $m$-dim subspace $\distr(q)$ of $T_qM$ and $\forall~q\in M$ we have
\bqn\llabel{Hor}\span{[X_1,[\ldots[X_{k-1},X_k]\ldots]](q)~|~X_i\in\mathrm{Vec}_H(M)}=T_qM
\eqn
where $\mathrm{Vec}_H(M)$ denotes the set of {\bf horizontal smooth vector fields} on $M$, i.e. $$\mathrm{Vec}_H(M)=\Pg{X\in\mathrm{Vec}(M)\ |\ X(p)\in\distr(p)~\ \forall~p\in M}.$$
\i $\g_q$ is a Riemannian metric on $\distr(q)$, that is smooth 
as function of $q$.
\ei
When $M$ is an orientable manifold, we say that the sub-Riemannian manifold is orientable.
\edeff

A Lipschitz continuous curve $\ga:[0,T]\to M$ is said to be \b{horizontal} if 
$\dot\ga(t)\in\distr(\ga(t))$ for almost every $t\in[0,T]$. Given an horizontal curve $\ga:[0,T]\to M$, the {\it length of $\ga$} is
\bqn
l(\ga)=\int_0^T \sqrt{ \g_{\ga(t)} (\dot \ga(t),\dot \ga(t))}~dt.
\eqnl{e-lunghezza}
The {\it distance} induced by the sub-Riemannian structure on $M$ is the 
function
\bqn
d(q_0,q_1)=\inf \{l(\ga)\mid \ga(0)=q_0,\ga(T)=q_1, \ga\ \mathrm{horizontal}\}.
\eqnl{e-dipoi}

The \hp\ of connectedness of M and the H\"ormander condition guarantee the finiteness and the continuity of $d(\cdot,\cdot)$ with respect to the topology of $M$ (Chow's Theorem, see for instance \cite{agra-book}). The function $d(\cdot,\cdot)$ is called the Carnot-Charateodory distance and gives to $M$ the structure of metric space (see \cite{bellaiche,gromov}).

Locally, the pair $(\distr,{\mathbf g})$ can be given by assigning a set of $m$ smooth vector fields spanning $\distr$ and that are orthonormal for ${\mathbf g}$, i.e.  
\bqn
\distr(q)=\span{X_1(q),\dots,X_m(q)}, ~~~\metr_q(X_i(q),X_j(q))=\delta_{ij}.
\eqnl{trivializable}
In this case, the set $\Pg{X_1,\ldots,X_m}$ is called a local {\bf orthonormal frame} for the sub-Riemannian structure. When  $(\distr,{\mathbf g})$ can be defined as in \r{trivializable} by $m$ vector fields defined globally, we say that the sub-Riemannian manifold is {\it trivializable}. 

%%%%%%%%%%%%%%%%%%%%%%%%%%%%%%%%%%%%%%%%%%%%%%%%%%%%%%%%%%%%%%%%%%%%

When the manifold is analytic and the orthonormal frame can be assigned through $m$ analytic vector fields, we say that the sub-Riemannian manifold is {\it analytic}.

We end this section with the definition of regular sub-Riemannian manifold.
\bdeff
Let $\distr$ be a distribution and define through the recursive formula
$$\distr_1:=\distr,~~~~~\distr_{n+1}:=\distr_n+[\distr_n,\distr].$$
%where 
%\bqn
%\distr_{n+1}(q_0)&:=&\distr_n(q_0)+[\distr_n,\distr](q_0)=\nn
%&=&\Pg{X_1(q_0)+[X_2,X_3](q_0)\ |\ X_1(q),X_2(q)\in\distr_n(q),\ X_3(q)\in\distr(q)~~\forall~q\in M}.
%\eqnn
The small flag of $\distr$ is the sequence 
$$\distr_1\subset\distr_2\subset\ldots\subset\distr_n\subset\ldots$$

A sub-Riemannian manifold is said to be {\bf regular} if for each $n=1,2,\ldots$ the dimension of $\distr_n(q_0)$ %=\Pg{f(q_0)\ |\ f(q)\in\distr_n(q)~\forall\ q\in M}
does not depend on the point $q_0\in M$.
\edeff

In this paper we always deal with sub-Riemannian manifolds that are orientable, analytic, trivializable and regular.

\subsection{Left-invariant sub-Riemannian manifolds}
\label{ss-leftmanifold}

In this section we present a natural sub-Riemannian structure that can be defined on Lie groups. All along the paper, we use the notation for Lie groups of matrices. For general Lie groups, by $gv$ with $g\in \G$ and $v\in \L$, we mean $(L_g)_*(v)$ where $L_g$ is the left-translation of the group.

\bdeff 
Let $\G$ be a Lie group with Lie algebra $\L$ and $\P\subseteq\L$ a subspace of $\L$ satisfying the {\bf Lie bracket generating condition} $$\Lie~\P:=\span{[\p_1,[\p_2,\ldots,[\p_{n-1},\p_n]]]\ |\ \p_i\in\P}=\L.$$
Endow $\P$ with a positive definite quadratic form $\Pa{.,.}$. Define a sub-Riemannian structure on $\G$ as follows:
\bi
\i the distribution is the left-invariant distribution $\distr(g):=g\P$;
\i the quadratic form $\metr$ on $\distr$ is given by $\metr_g(v_1,v_2):=\Pa{g^{-1}v_1,g^{-1}v_2}$.
\ei
In this case we say that $(\G, \distr, \metr)$ is a left-invariant sub-Riemannian manifold.
\llabel{d-lieg-leftinv}
\edeff
\brem Observe that all left-invariant manifolds $(\G, \distr, \metr)$ are regular.
\erem

In the following we define a left-invariant sub-Riemannian manifold choosing a set of $m$ vectors $\Pg{\p_1,\ldots,\p_m}$ that are an orthonormal basis for the subspace $\P\subseteq\L$ with respect to the metric defined in Definition \ref{d-lieg-leftinv}, i.e. $\P=\span{\p_1,\ldots,\p_m}$ and $\Pa{\p_i,\p_j}=\de_{ij}$. We thus have $\distr(g)=g\P=\span{g\p_1,\ldots,g\p_m}$ and $\g_g(g\p_i,g\p_j)=\de_{ij}$. Hence, every left-invariant sub-Riemannian manifold is trivializable.

\subsubsection{The intrinsic hypoelliptic Laplacian}
\label{sez-intr}
\newcommand{\muh}{\mu_\SR}
\newcommand{\divh}{\mathrm{div}_\SR}

In this section, we recall the definition of intrinsic hypoelliptic Laplacian 
given in \cite{nostro-kern} and based on the Popp volume form in sub-Riemannian geometry presented in \cite{montgomery}.

Let $(M,\distr,\g)$ be a $(n,m)$-sub-Riemannian manifold and $\{X_1,\ldots X_m\}$ a local orthonormal frame. The operator obtained by the sum of squares of these vector fields is not a good definition of hypoelliptic Laplacian, since it depends on the choice of the orthonormal frame (see for instance \cite{nostro-kern}).

In sub-Riemannian geometry an invariant definition of  hypoelliptic Laplacian is obtained by computing the divergence of the horizontal gradient, like the Laplace-Beltrami operator in Riemannian geometry.

\bdeff Let  $(M,\distr,\g)$ be an orientable  regular sub-Riemannian manifold. We define the intrinsic hypoelliptic Laplacian as
$\Dh\phi:=\divh\gradh\phi$, where 

\bi
\i the horizontal gradient is the unique operator $\gradh$ from $\mathcal{C}^\infty(M)$ to ${{\mathrm{Vec}}_H(M)}$ 
satisfying $\g_q(\gradh\phi(q),v)=d\phi_q (v)~~~~\forall~q\in M,~v\in \distr(q)$. (In coordinates if  $\{X_1,\ldots X_m \}$ is a local orthonormal frame for $(M,\distr,\metr)$, then $\gradh \phi=\sum_{i=1}^m \Pt{L_{X_i}\phi}X_i$.)
\i the divergence of a vector field $X$ is the unique function satisfying   
$\div X \muh=L_X \muh$ where $\muh$ is the Popp volume form.
\ei
\edeff
The construction of the Popp volume form is not totally trivial and we address the reader to \cite{montgomery} or \cite{nostro-kern} for details. We just recall that the Popp volume form coincide with the Lebesgue measure in a special system of coordinate related to the nilpotent approximation. In sub-Riemannian geometry one can also define other intrinsic volume forms, like the Hausdorff or the spherical Hausdorff volume. However, at the moment, the Popp volume form is the only one  known to be smooth in general. However for left-invariant sub-Riemannian manifolds  all these measures are proportional to the left Haar measure. 

The hypoellipticity of $\Dh$ (i.e. given $U\subset M$ and $\fz{\phi}{U}{\R}$ such that $\Dh\phi\in\cont^\infty$, then $\phi$ is $\cont^\infty$) follows from the H\"ormander Theorem (see \cite{hormander}).

In this paper we are interested only to nilpotent Lie groups. The next proposition says that for all unimodular Lie groups, i.e. for groups such that the left and right Haar measure coincides (and in particular for real connected nilpotent groups) the intrinsic hypoelliptic Laplacian is the sum of squares.

\bp
Let $(\G,\distr,\metr)$ be a left-invariant sub-Riemannian manifold generated by the orthonormal basis $\Pg{\p_1,\ldots,\p_m}\subset\l$. 
If $G$ is unimodular then $\Dh\phi=\sum_{i=1}^m \Pt{L_{X_i}^2\phi}$
where $L_{X_i}$ is the Lie derivative w.r.t. the field $X_i=g\p_i$.
\ep

\subsection{Computation of the hypoelliptic heat kernel via the Generalized Fourier Transform}
\label{s-genfour}

\newcommand{\Rep}{\Repbase^\lam}
\newcommand{\dRep}{\dRepbase^\lam}
\newcommand{\domR}{\domRbase^\lam}
\newcommand{\fzR}{\fzRbase^\lam_n}
\newcommand{\eigR}{\eigRbase^\lam_n}
\newcommand{\DHL}{\DHLbase^\lam}
\renewcommand{\G}{G}

In this section we describe the method, developed in \cite{nostro-kern}, for the computation of the hypoellitpic heat kernel for left-invariant sub-Riemannian structures on unimodular Lie groups.

The method is based upon the \GFT, that permits to disintegrate a function from a Lie  group $\G$ to $\R$ on its components on (the class of) non-equivalent unitary irreducible representations of $\G$. For proofs and more details, see \cite{nostro-kern}.

\subsubsection{The Generalized Fourier Transform}

Let $f\in L^1(\R,\R)$: its Fourier transform is defined by the formula
$$
\hat f(\lam)=\int_\R f(x) {e^{-i x \lam}} dx.
$$
If $f\in L^1(\R,\R) \cap L^2(\R,\R)$ then $\hat f\in L^2(\R,\R)$ and one has
$$
\int_\R|f(x)|^2dx=\int_\R|\hat f(\lam)|^2  \frac{d\lam}{2\pi},
$$
called Parseval or Plancherel equation. By density of $L^1(\R,\R) \cap L^2(\R,\R)$ in $L^2(\R,\R)$, this equation expresses the fact that the Fourier transform is an isometry between $L^2(\R,\R)$ and itself. Moreover, the following inversion formula holds:
$$
f(x)=\int_\R \hat f(\lam) e^{i x \lam} \frac{d\lam}{{2 \pi}},
$$
where the equality is intended in the $L^2$ sense. 
It has been known from more than 50 years that the Fourier transform generalizes to a wide class of locally compact groups (see for instance \cite{chirichian,duflo,hewittI,hewittII,kirillov-el,taylor-GFT}).  Next we briefly present this generalization for groups satisfying the following \hp:
\bd

\i  \Hp\ $\G$ is a unimodular Lie group of Type I.
\ed
For the definition of groups of Type I see \cite{dixmier1}. For our purposes it is sufficient to recall that all groups treated in this paper (i.e. \gruppi) are of Type I. Actually, all the real connected nilpotent Lie groups are of Type I \cite{dixmier,harish}.
In the following, the $L^p$ spaces $L^p(\G,\C)$ are intended with respect to the Haar measure $\mu:=\mu_L=\mu_R$.
 
Let $\G$ be a Lie group satisfying \Hp\ and $\hat \G$ be the dual\footnote{\label{nota-dual} In this paper, by the dual of the group, we mean  the support of the Plancherel measure on the set of non-equivalent unitary irreducible representations of $G$; we thus ignore the singular representations.} of the group $\G$, i.e. the set of all equivalence classes of unitary irreducible representations of $\G$. Let $\lam\in\hat \G$: in the following we indicate by $\Rep$ a choice of an irreducible representation in the class $\lam$. By definition, $\Rep$ is a map that to an element of $\G$ associates a unitary operator acting on a complex separable  Hilbert space ${\H^\lam}$:
\mmfunz{\Rep}{\G}{U(\domR)}{g}{\Rep(g).}

The index $\lam$ for $\H^\lam$ indicates that in general the Hilbert space can vary with $\lam$.
 
\bdeff
Let $\G$ be a Lie group satisfying \Hp, and $f\in L^1(\G,\C)$. The generalized (or noncommutative) Fourier transform (\GFT) of $f$ is the map (indicated in the following as $\hat f$ or $\FF(f)$)
that to each element of $\hat \G$ associates the linear operator on $\domR$:
 \bqn
 \hat f(\lam):=\FF(f):=\int_\G    f(g)\Rep(g^{-1})d\mu.
 \eqn
  \edeff
\noindent
Notice that since $f$ is integrable and  $\Rep$ unitary, then $\hat f(\lam)$ is a bounded operator.

\brem
$\hat f$ can be seen as an operator from $\stackrel{\oplus}{\int_{\hat \G}}
\H^\lam$ to itself. We also use the notation  $\hat f=\stackrel{\oplus}{\int_{\hat \G}}\hat f(\lam)$
\erem

In general $\hat \G$ is not a group and its structure can be quite complicated.  In the case in which $\G$ is abelian then $\hat \G$ is a group; if $\G$ is nilpotent (as in our cases) then $\hat \G$ has the structure of $\R^n$ for some $n$.

Under the \hp\ \Hp\ one can define on $\hat \G$ a positive measure $dP(\lam)$   (called the Plancherel measure)  such that for every $f\in L^1(\G,\C)\cap L^2(\G,\C)$ one has 
$$
\int_\G |f(g)|^2 \mu(g)=\int_{\hat \G}Tr(\hat f(\lam)\circ \hat f (\lam)^\ast ) dP(\lam).
$$
By density of $L^1(\G,\C)\cap L^2(\G,\C)$ in $L^2(\G,\C)$, this formula expresses the fact that the \GFT\ is an isometry between $L^2(\G,\C)$ and 
$\stackrel{\oplus}{\int_{\hat \G}} \HS^\lam$, the set  of Hilbert-Schmidt operators with respect to the Plancherel measure.
Moreover, it is obvious that:
\bp
Let $\G$ be a Lie group satisfying \Hp and $f\in L^1(\G,\C)\cap L^2(\G,\C)$. We have, for each $g\in \G$
\bqn
f(g)=\int_{\hat \G}Tr(\hat f (\lam)\circ   \Rep(g) ) dP(\lam).
\eqn  where the equality is intended in the $L^2$ sense. 
\ep
\noindent
{It is immediate to verify that, given two functions $f_1,f_2\in  L^1(\G,\C)$ and defining their convolution as
\bqn
(f_1\ast f_2)(g)=\int_\G  f_1(h)f_2(h^{-1}g) dh,
\eqn
then the GFT maps  the convolution into non-commutative product:  
\bqn 
\FF( f_1\ast f_2) (\lam)=\hat{f}_2 (\lam)\hat{f}_1(\lam).
\label{conv}
\eqn
Another important property is that if $\delta_{\Id}(g)$ is the Dirac function at the identity over $\G$, then 
\bqn
\hat{\delta}_{\Id}(\lam)=\Id_{H^\lam}.
\eqn
}
In the following, a key role is played by  the infinitesimal version of the representation $\Rep$, that is the map 
\bqn
\label{diff-rep}
d\Rep:X\mapsto d\Rep(X):=\left.\frac{d}{dt}\right|_{t=0}\Rep(e^{t p}),
\eqn
where $X=gp$, ($p\in\l$, $g\in \G$) is a left-invariant vector field over $\G$. By Stone theorem (see for instance \cite[p. 6]{taylor-GFT})  $d\Rep(X)$ is a (possibly unbounded) skew-adjoint operator on 
%the tangent of $\H_\lam$ that we identify with 
$\H^\lam$. 
We have the following:
\bp
\label{p-uffa}
Let $\G$ be a Lie group satisfying \Hp\ and $X$ be a left-invariant vector field over $\G$. The GFT of $X$, i.e. $\hat X=\FF L_X\FF^{-1}$ splits into the Hilbert sum of operators $\hat X^\lam$, \gr{each one of them acting} on the set $\HS^\lam$ of Hilbert-Schmidt operators over $\H^\lam$:
  \bqn
 \hat X=\stackrel{\oplus}{\int_{\hat \G}}\hat X^\lam.
 \eqnn
Moreover,
  \bqn
 \hat X^\lam\Xi=d\Rep(X)\circ\Xi,~~\mbox{ for every }\Xi\in \HS^\lam,
 \eqn
 i.e. the GFT of a left-invariant vector field acts as a left-translation over $\HS^\lam$.
\ep
\noindent

\brem
\label{r-accalambda}
From the fact that the GFT of a left-invariant vector field acts as a left-translation, it follows that $\hat X^\lam$ can be interpreted as an operator over $\H^\lam$.
\erem

\subsubsection{Computation of the kernel of the hypoelliptic heat equation}
\label{s:method}
\newcommand{\du}{L}

In this section we provide a general method  to compute the kernel of the hypoelliptic heat equation on a left-invariant sub-Riemannian manifold $(\G,\distr,\metr)$ such that $\G$ satisfies the assumption \Hp. 

We begin by recalling some existence results (for the semigroup of evolution and for the corresponding kernel) in the case of the  sum of squares. We recall that for all the examples treated in this paper the invariant hypoelliptic Laplacian is the sum of squares. 

Let $\G$ be a unimodular Lie group and $(\G,\distr,\metr)$ a left-invariant sub-Riemannian manifold generated by the orthonormal basis $\Pg{\p_1,\ldots,\p_m}$, and consider the hypoelliptic heat equation
\bqn
\partial_t\phi(t,g)=\Dh\phi(t,g).
\llabel{eq-hypoQ}
\eqn
Since $\G$ is unimodular, then $\Dh=L_{X_1}^2+\ldots +L_{X_m}^2$, where $L_{X_i}$ is the Lie derivative w.r.t. the vector field $X_i:=g\p_i$ $(i=1,\ldots, m)$.
Following Varopoulos \cite[pp. 20-21, 106]{varopoulos}, since $\Dh$ is a sum of squares, then it is a symmetric operator that we identify with its Friedrichs (self-adjoint) extension, that is the infinitesimal generator of a (Markov) semigroup $e^{t \Dh }$. Thanks to the left-invariance of $X_i$ (with $i=1,\ldots, m)$, $e^{t \Dh }$ admits a a right-convolution kernel $p_t(.)$, i.e. 
\bqn
 e^{t \Dh } \phi_0(g)=\phi_0\ast p_t (g)=\int_\G  \phi_0(h)p_t(h^{-1}g) \mu(h)
\eqn
is the solution for $t>0$ to $\r{eq-hypoQ}$ with initial condition $\phi(0,g)=\phi_0(g)\in L^1(\G,\R)$ with respect to the Haar measure.

Since the operator $\partial_t-\Dh$ is hypoelliptic, then the kernel is a $\cont^\infty$ function of $(t,g)\in \R^+\times \G$. 
Notice that $p_t(g)=e^{t \Dh }\de_\Id(g)$.

The main results of the paper are based on the following key fact.
\bt
\label{t-main}
Let $\G$ be a Lie group satisfying \Hp\ and $(\G,\distr,\metr)$  a left-invariant sub-Riemannian manifold generated by the orthonormal basis $\Pg{\p_1,\ldots,\p_m}$. Let $\Dh=L_{X_1}^2+\ldots +L_{X_m}^2$ be the intrinsic hypoelliptic Laplacian where $L_{X_i}$ is the Lie derivative w.r.t. the vector field $X_i:=g\p_i$. 

Let $\Pg{\Rep}_{\lam\in \hat \G}$ be the set of all non-equivalent classes of irreducible representations of the group $\G$, each acting on an Hilbert space $\domR$, and $dP(\lam)$ be the Plancherel measure on the dual space $\hat \G$. We have the following:
\be[{(\bf i)}]
\i the GFT of $\Dh$ splits into the Hilbert sum of operators $\DHL$, \gr{each one of which leaves}  ${\H^\lam}$ invariant:
\bqn
\llabel{oplus}
\hat{\Delta}_\SR=\FF\Dh\FF^{-1}=\stackrel{\oplus}{\int_{\hat \G}}\DHL dP(\lam),
\mbox{~~ where ~~}
\DHL=\sum_{i=1}^m \left(\hat{X}_i^\lam\right)^2.
\eqn
\i The operator $\DHL$ is self-adjoint and it is the infinitesimal generator of a contraction semi-group $e^{t\DHL}$ over $\HS^\lam$, i.e. $e^{t\DHL}  \Xi_0^\lam$ is the solution for $t>0$ to the operator equation $\partial_t  \Xi^\lam(t)=\DHL \Xi^\lam(t)$ in $\HS^\lam$, with initial condition  $\Xi^\lam(0)= \Xi^\lam_0 $.
\i The hypoelliptic heat kernel is
\bqn
\llabel{formula-1}
p_t(g)=\int_{\hat \G}Tr\left(e^{t\DHL}\Rep(g)\right)dP(\lam),~~t>0.
\eqn
\ee
\et
\noindent
\brem
As a consequence of Remark \ref{r-accalambda}, it follows that $\DHL$ and $e^{t\DHL}$ can be considered as operators on $\H^\lam$.

\erem

%\brem
%\label{ispazi}
%Notice that $\DHL$ is an operator on $U(\H_\lam)$ i.e.  $\DHL\hat f(\lam)$ (where $f\in L^1(G,\C)$)) is an operator on $\H_\lam$. Hence  $\DHL $ can be though as a self-adjoint operator acting on $\H_\lam$ on vectors of the type  $\Xi^\lam=\hat f(\lam)\psi^\lam$, where $\psi^\lam\in\H_\lam$. The same remark holds for $e^{t\DHL}$.
%In the following, when we speak of eigenvalues and eigenvectors of $\DHL$, we think to  $\DHL$ as an operator on  $\H_\lam$. 
%\erem

%\gr{In the case when each $\DHL$ has discrete spectrum, the following corollary gives an explicit formula for the hypoelliptic heat kernel in terms of its eigenvalues and eigenvectors.}

%\bc Under the hypotheses of Theorem \ref{t-main}, if in addition \gr{we have that} for every $\lam$, $\DHL$ (considered as an operator over $\H^\lam$) has discrete spectrum and  $\Pg{\fzR}$ is a complete set of eigenfunctions of norm one with the corresponding set of eigenvalues $\Pg{\eigR}$, then
%\bqn
%p_t(g)=\int_{\hat \G}\left( \sum_{n} e^{\eigR t} \Pa{\fzR ,\Rep(g) \fzR} \right)d P(\lam)
%\eqnl{e-fundsol-general}
%where $\Pa{.,.}$ is the scalar product in $\domR$.
%\ec

The following corollary gives a useful formula for the hypoelliptic heat kernel in the case in which for all $\lam\in\hat{\G}$ each operator $e^{t\DHL}$ admits a convolution kernel $Q_t^\lam(.,.)$. Below by $\fzRbase^\lam$, we \gr{intend} an element of $\H^\lam$.

\bc
\label{c-2} Under the hypotheses of Theorem \ref{t-main}, if for all $\lam\in\hat \G$ we have $\H^\lam=L^2(X^\lam,d\th^\lam)$ for some measure space $(X^\lam,d\th^\lam)$ and 
$$\Pq{e^{t\DHL} \fzRbase^\lam}(\th)=\int_{X^\lam}  \fzRbase^\lam (\bar \th) Q_t^\lam(\th,\bar\th)\,d\bar \th,$$ then
%%%%%%%%%%%%%%%%%%%%%%%%%%%
$$
p_t(g)=\int_{\hat \G} \int_{X^\lam}\left.\Rep(g)Q_t^\lam(\th,\bar\th)\right|_{\th=\bar \th}              \,d\bar \th\,
d P(\lam),
$$
where in the last formula $\Rep(g)$ acts on $Q_t^\lam(\th,\bar\th)$ as a function of $\th$.
\ec
%\bproof  
%From \r{formula-1}, we have
%$$
%p_t(g)=\int_{\hat \G}Tr\left(e^{t\DHL}\Rep(g)\right)dP(\lam)=\int_{\hat \G} Tr\left(\Rep(g)e^{t\DHL}\right)dP(\lam).
%$$
%We have to compute the trace of the operator
%\newcommand{\oprusso}{\Theta}
%\bqn
%\label{kernellino}
%\oprusso=\Rep(g)e^{t\DHL}:\fzRbase^\lam(\th)\mapsto \Rep(g)e^{t\DHL}\fzRbase^\lam(\th)&=&
%\Rep(g)\int_{X^\lam}  \fzRbase^\lam(\bar \th) Q_t^\lam(\th,\bar\th)\,d\bar \th=\\
%&=&\int_{X^\lam}  K(\th,\bar\th) \fzRbase^\lam(\bar \th)d\bar \th
%\eqnn
%where  $K(\th,\bar\th)= \Rep(g)Q_t^\lam(\th,\bar\th)$ is a function of $\th,\bar\th$ and $\Rep(g)$ acts on $Q_t^\lam(\th,\bar\th)$ as a function of $\th$. The trace of $\oprusso$ is $\int_X K(\bar \th,\bar \th)d\bar\th$ and the conclusion follows.
%\eproof

\renewcommand{\G}{\mathfrak{g}}

\section{Hypoelliptic heat kernels on \gruppi}

\label{s-gruppi}
In this section we describe the groups \gruppi\ and we provide their matrix and Euclidean presentations. We define left-invariant sub-Riemannian structures on them and the corresponding hypoelliptic Laplacian.

We then provide representations of the groups and compute the GFT of the hypoelliptic Laplacian. We apply the method presented in Section \ref{s:method} to compute the fundamental solution of the hypoelliptic heat equation.

\subsection{Definitions of \gruppi}

In our paper we deal with two 3-step Lie groups. The first one is the nilpotent group $\GG_4$ with growth vector $(2,3,4)$. Its Lie algebra is $\L_4=\mathrm{span}\Pg{\l_1,\l_2,\l_3,\l_4},$ whose  generators satisfy
\bqn
&&\Pq{\l_1,\l_2}=\l_3,\ \Pq{\l_1,\l_3}=\l_4,\ \Pq{\l_1,\l_4}=\Pq{\l_2,\l_3}=\Pq{\l_2,\l_4}=\Pq{\l_3,\l_4}=0.
\eqnn 

The second one is the free nilpotent group $\GG_5$ with growth vector $(2,3,5)$. Its Lie algebra is $\L_5=\mathrm{span}\Pg{\l_1,\l_2,\l_3,\l_4,\l_5},$ whose generators satisfy
\bqn
&&\Pq{\l_1,\l_2}=\l_3,\ \Pq{\l_1,\l_3}=\l_4,\ \Pq{\l_2,\l_3}=\l_5,\nn
&&\Pq{\l_1,\l_4}=\Pq{\l_1,\l_5}=\Pq{\l_2,\l_4}=\Pq{\l_2,\l_5}=\Pq{\l_3,\l_4}=
\Pq{\l_3,\l_5}=\Pq{\l_4,\l_5}=0.
\eqnn

\ffoot{We recall here what means that the groups defined above are 3-step nilpotent. Consider a Lie algebra $\L$  and define the {\it lower central series} recursively as follows:
$$\L^1:=\L\qquad \L^{n+1}:=\Pq{\L^n,\L}=\Pg{\Pq{l,m}\ \mid\ l\in\L^n,\,m\in\L}.$$
As a direct consequence, we have that $\L^{n}\supset\L^{n+1}$. We have the following standard definition.
\bdeff
A Lie algebra $\L$ is nilpotent if it satisfies $\L^n=0$ for some $n\in\N$.

A nilpotent Lie algebra $\L$ is $n$-step if $\L^{n+1}=0$ and $\L^n\neq 0$.

A Lie group is ($n$-step) nilpotent if its Lie algebra is ($n$-step) nilpotent.
\edeff}
Both \gruppi\ are 3-step nilpotent, as a direct consequence of the definition.

\subsection{Hypoelliptic heat kernel on $\GG_4$}

In this section we first give the matrix and Euclidean presentations of the Lie group $\GG_4$. We then define a sub-Riemannian structure on it. We give explicitly the representations of the group, that we use at the end to compute the hypoelliptic kernel in terms of the kernel of the quartic oscillator.

We start with the Lie algebra $\L_4$, that can be presented as the follow matrix space
$$\L_4\simeq\Pg{\Mat{cccc}{
0&-a_1& 0& a_4\\
0&0  & -a_1& a_3\\
0&0  & 0& a_2\\
0&0  & 0& 0}\ \mid\ a_i\in\R}.$$
We present each $\l_i$ as the matrix with $a_j=\de_{ij}$. It is straightforward to prove that these matrices satisfy the commutation rules for $\L_4$, where the  bracket operation is the standard $\Pq{A,B}:=BA-AB$.

A matrix presentation of the group $\GG_4$ is thus the matrix exponential of $\L_4$:
\bqn\GG_4&\simeq&\Pg{\exp{\Mat{cccc}{
0&-a_1& 0& a_4\\
0&0  & -a_1& a_3\\
0&0  & 0& a_2\\
0&0  & 0& 0}}\ \mid\ a_i\in\R}=\Pg{\Mat{cccc}{
1&-x_1& \frac{x_1^2}{2}& x_4\\
0&1  & -x_1& x_3\\
0&0  & 1& x_2\\
0&0  & 0& 1}\ \mid\ x_i\in\R},
\eqnn
with
\bqn
x_1&=&a_1,\qquad x_2=a_2,\qquad x_3=a_3-\frac{a_1a_2}{2},\quad x_4=a_4+\frac{a_1^2a_2}{6}-\frac{a_1a_3}{2}.
\eqnn
%Thus, $\GG_4$ is a group with respect to the standard matrix product.

We now define the isomorphism $\Pi_4$ between $\GG_4$ and $\R_4$ given by $$\Pi_4\Pt{\Mat{cccc}{
1&-x_1& \frac{x_1^2}{2}& x_4\\
0&1  & -x_1& x_3\\
0&0  & 1& x_2\\
0&0  & 0& 1}}=(x_1,x_2,x_3,x_4).$$
This isomorphism is a group isomorphism when $\R^4$ is endowed with the following product (see \cite[p. 330]{dixmier3}): \bqn(x_1,x_2,x_3,x_4)\cdot(y_1,y_2,y_3,y_4)&:=&
\Pt{x_1 + y_1, x_2 + y_2, x_3+y_3-x_1 y_2,
x_4+y_4+\frac12 x_1^2 y_2 - x_1 y_3}
\eqnn
%We thus have that $(\R^4,\cdot)$ is a Lie group. 
The isomorphism $\Pi_4$ induces an isomorphism of tangent spaces $T_g \G \simeq T_{\Pi_4(g)} \R^4$, that is explicitly $g \l_i\simeq X_i$, with $X_i$ given by 
\bqn
X_1(x)&=&\der{x_1},\qquad
X_2(x)=\der{x_2} -x_1\der{x_3} + \frac{x_1^2}{2}\der{x_4},\label{e:Xi-G4}\\
X_3(x)&=&\der{x_3} -x_1\der{x_4},\qquad
X_4(x)=\der{x_4},
\eqnn
where $x=(x_1,x_2,x_3,x_4)$.

\subsubsection{Left-invariant sub-Riemannian structure on $\GG_4$}

We endow $\GG_4$ with a left-invariant sub-Riemannian structure as presented in Section \ref{ss-leftmanifold}. We define the sub-Riemannian manifold $(\GG_4,\distr,\metr)$ where $\distr(g)=g\p$ with $\p=\span{\l_1,\l_2}$ and $\metr_g(g\l_i,g\l_j)=\de_{ij}$ with $i,j=1$ or $2$. %This structure is indeed unique. See  \cite{yuri4}.

Since $\GG_4$ is nilpotent, then it is unimodular, thus the intrinsic hypoelliptic Laplacian $\Delta_H$ is the sum of squares (see \cite[Proposition 17]{nostro-kern}). In terms of the Euclidean presentation of $\GG_4$, the hypoelliptic Laplacian is thus $\Delta_H=X_1^2+X_2^2$, with the $X_i$ given by \r{e:Xi-G4}.

We thus want to find the fundamental solution for the following heat equation:
\bqn
\partial_t \phi(t,x)=\Dh \phi(t,x).
\eqnl{e:hypo-G4}

\subsubsection{Representations of $\GG_4$}

\renewcommand{\Rep}{\Repbase^{\lam,\mu}}
\renewcommand{\dRep}[1]{\dRepbase_{#1}^{\lam,\mu}}
\renewcommand{\domR}{\domRbase}
\renewcommand{\fzR}{\fzRbase}
\renewcommand{\DHL}{\DHLbase^{\lam,\mu}}

We now recall the representations of the group $\GG_4$, as computed by Dixmier in \cite[p. 333]{dixmier3}. As stated before, we may consider only representations on the support of the Plancherel measure.
\bp
The dual space of $\GG_4$ is $\hat{G}=\Pg{\Rep\ |\ \lam\neq 0,\mu\in\R}$, where
\mmfunz{\Rep(x_1,x_2,x_3,x_4)}{\domR}{\domR}{\fzR(\th)}{\exp{i\Pt{-\frac{\mu}{2\lam}x_2+\lam x_4 - \lam x_3 \th + \frac{\lam}{2} x_2 \th^2}}\fzR(\th+x_1)}
whose domain is $\domR=L^2(\R,\C)$, endowed with the standard product $<\fzR_1,\fzR_2>:=\int_\R \fzR_1(\th)\overline{\fzR_2}(\th)\,d\th$ where $d\th$ is the Lebesgue measure.

The Plancherel measure on $\hat{G}$ is $dP(\lam,\mu)=d\lam d\mu$, i.e. the Lebesgue measure on $\R^2$.
\ep
\brem
Notice that in this case the domain $\domR$ of the representation $\Rep$ does not depend on $\lam,\mu$.
\erem

\subsubsection{The kernel of the hypoelliptic heat equation}

Consider the representation $\Rep$ of $\GG_4$ and let $\dRep{i}$ be the corresponding representations of the differential operators $L_{X_i}$ with $i=1,2$. Recall that $\dRep{i}$  are operators on $\domR$. Again from \cite[p. 333]{dixmier3}, or by explicit computation, we have
\bqn
\Pq{\dRep{1} \fzRbase}(\th)=\frac{d}{d\th}\fzRbase(\th),\quad
[\dRep{2}\fzRbase](\th)=
\Pt{-\frac{i}{2}\frac{\mu}{\lam}+\frac{i}2 \lam \th^2}\fzRbase(\th),\eqnn
thus
\bqn
\Pq{\DHL\fzRbase}(\th)=\Pt{\frac{d^2}{d\th^2}-\frac14\Pt{\lam\th^2-\frac{\mu}{\lam}}^2}\fzRbase(\th).\eqnn
The GFT of the hypoelliptic heat equation is thus
\bqn
\partial_t\fzRbase= \Pt{\frac{d^2}{d\th^2}-\frac14\Pt{\lam\th^2-\frac{\mu}{\lam}}^2}\fzRbase(\th).\eqnl{e:GFT-G4}
We rewrite it as
\bqn
\partial_t\fzRbase= \Pt{\frac{d^2}{d\th^2}-\Pt{\al\th^2+\beta}^2}\fzRbase(\th),\eqnl{e:GFT-G4-trasf}
with $\al=\frac\lam2$, $\beta=-\frac{\mu}{2\lam}$.

\newcommand{\sol}{\Psi}

The operator $\frac{d^2}{d\th^2}-(\al\th^2+\beta)^2$ is the Laplacian with quartic potential, see e.g. \cite{skala}. As already stated, no general explicit solutions are known for this equation. We call 
\bqn\sol_t \Pt{\th,\bar\th;\al,\beta}
\eqnn
the solution of
\bqn
\begin{cases}
\partial_t \psi(t,\th)=\Pt{\frac{d^2}{d\th^2}-\Pt{\al\th^2+\beta}^2}\psi(t,\th),\\
\psi(0,\th)=\de_{\bar\th},
\end{cases}
\eqnn
 i.e. the solution of \r{e:GFT-G4-trasf} evaluated in $\th$ at time $t$, with initial data $\de_{\bar\th}$ and parameters $\al,\beta$.

Applying Corollary \ref{c-2} and after straightforward computations, one gets  the kernel of the hypoelliptic heat equation on the group $\GG_4$:
\bqn
p_t(x_1,x_2,x_3,x_4)=\int_{\R\backslash\Pg{0}} d\lam \int_\R d\mu \int_\R d \th\,
e^{i\Pt{-\frac{\mu}{2\lam}x_2+\lam x_4 - \lam x_3 \th + \frac{\lam}{2} x_2 \th^2}}\sol_t\Pt{\th+x_1,\th;\frac\lam2,-\frac{\mu}{2\lam}}.
\eqnl{eq-G4-heat-exp}

%%%%%%%%%%%%%%%%%%%%%%%%%%%%%%%%%% G5 %%%%%%%%%%%%%%%%%%%%%%%%%%%%%%%%%%%%%

\subsection{Hypoelliptic heat kernel on $\GG_5$}

%In this section we give the matrix and Euclidean presentations of the Lie group $\GG_5$. We then define a sub-Riemannian structure on it. We give explicitly the representations of the group, that we use at the end to compute the hypoelliptic kernel in terms of the kernel of the quartic oscillator.

The Lie algebra $\L_5$ of the group $\GG_5$ can be presented as the following matrix space
$$\L_5\simeq\Pg{\Mat{cc}{\b{M_1}(a_1,a_2,a_3,a_4) & \b{0}_{4\times 4}\\
\b{0}_{4\times 4} & \b{M_2}(a_1,a_2,a_3,a_5)}\ |\ 
 a_i\in\R},$$
where
$$
\b{M_1}(a_1,a_2,a_3,a_4)=\Mat{cccc}{
0&-a_1& 0& a_4\\
0&0  & -a_1& a_3\\
0&0  & 0& a_2\\
0&0  & 0& 0},\qquad
 \b{M_2}(a_1,a_2,a_3,a_5)=\Mat{cccc}{
0&a_2& 0& a_5\\
0&0  & a_2& -a_3\\
0&0  & 0& -a_1\\
0&0  & 0& 0}.$$

We present each $\l_i$ as the matrix with $a_j=\de_{ij}$. It is straightforward to prove that these matrices satisfy the commutation rules for $\L_5$, where the  bracket operation is the standard $\Pq{A,B}:=BA-AB$.

A matrix presentation of the group $\GG_5$ is thus the matrix exponential of $\L_5$:
\bqn\GG_5&\simeq&\Pg{\exp{\Mat{cc}{\b{M_1}(a_1,a_2,a_3,a_4) & \b{0}_{4\times 4}\\
\b{0}_{4\times 4} & \b{M_2}(a_1,a_2,a_3,a_5)}}\ |\ a_i\in\R}=\nn
&=&\Pg{\Mat{cc}{\exp{\b{M_1}(a_1,a_2,a_3,a_4)} & \b{0}_{4\times 4}\\
\b{0}_{4\times 4} & \exp{\b{M_2}(a_1,a_2,a_3,a_5)}}\ |\ a_i\in\R}=\nn
&=&\Pg{
\Mat{cc}{\b{N_1}(x_1,x_2,x_3,x_4) & \b{0}_{4\times 4}\\
\b{0}_{4\times 4} & \b{N_2}(x_1,x_2,x_3,x_5)}\ |\ x_i\in\R}
\eqnn
with
\bqn
\b{N_1}(x_1,x_2,x_3,x_4)=\Mat{cccc}{
1&-x_1& \frac{x_1^2}{2}& x_4\\
0&1  & -x_1& x_3\\
0&0  & 1& x_2\\
0&0  & 0& 1},\quad
 \b{N_2}(x_1,x_2,x_3,x_5)=\Mat{cccc}{
1&x_2& \frac{x_2^2}{2}& x_5-\frac{x_1 x_2^2}2\\
0&1  & x_2& -x_3-x_1 x_2\\
0&0  & 1& -x_1\\
0&0  & 0& 1},\eqnn
\bqn
x_1&=&a_1,\qquad x_2=a_2,\qquad x_3=a_3-\frac{a_1a_2}{2},\quad x_4=a_4+\frac{a_1^2a_2}{6}-\frac{a_1a_3}{2},\quad x_5=a_5+\frac{a_1a_2^2}{6}-\frac{a_2a_3}{2}.
\eqnn
%Thus, $\GG_5$ is a group with respect to the standard matrix product.

We now define the isomorphism $\Pi_5$ between $\GG_5$ and $\R_5$ given by $$\Pi_4\Pt{\Mat{cc}{\b{N_1}(x_1,x_2,x_3,x_4) & \b{0}_{4\times 4}\\
\b{0}_{4\times 4} & \b{N_2}(x_1,x_2,x_3,x_5)}}=(x_1,x_2,x_3,x_4,x_5).$$
This isomorphism is a group isomorphism when $\R^5$ is endowed with the following product (see \cite[p. 331]{dixmier3}): \bqn(x_1,x_2,x_3,x_4,x_5)\cdot(y_1,y_2,y_3,y_4,y_5)&:=&\eqnn
\bqn
\Pt{x_1 + y_1, x_2 + y_2, x_3+y_3-x_1 y_2,
x_4+y_4+\frac12 x_1^2 y_2 - x_1 y_3,
x_5+y_5+\frac12 x_1 y_2^2 - x_2 y_3+ x_1 x_2 y_2}.
\eqnn
%We thus have that $(\R^4,\cdot)$ is a Lie group. 
The isomorphism $\Pi_5$ induces an isomorphism of tangent spaces $T_g \G \simeq T_{\Pi_5(g)} \R^5$, that is explicitly $g \l_i\simeq X_i$, with $X_i$ given by 
\bqn
X_1(x)&=&\der{x_1},\qquad
X_2(x)=\der{x_2} -x_1\der{x_3} + \frac{x_1^2}{2}\der{x_4}+x_1 x_2 \der{x_5},\label{e:Xi-G5}\\
X_3(x)&=&\der{x_3} -x_1\der{x_4} -x_2\der{x_5},\qquad
X_4(x)=\der{x_4},\qquad X_5(x)=\der{x_5},
\eqnn
where $x=(x_1,x_2,x_3,x_4,x_5)$.

%\subsubsection{Left-invariant sub-Riemannian structure on $\GG_5$}

We endow $\GG_5$ with a left-invariant sub-Riemannian structure as presented in Section \ref{ss-leftmanifold}, where $\p=\span{\l_1,\l_2}$ and $\metr_g(g\l_i,g\l_j)=\de_{ij}$ with $i,j=1$ or $2$.% This structure is indeed unique. See  \cite{yuri5}.

%Since $\GG_5$ is nilpotent, then it is unimodular, thus the intrinsic hypoelliptic Laplacian $\Delta_H$ is the sum of squares (see \cite[Proposition 17]{nostro-kern}). In terms of the Euclidean presentation of $\GG_5$, the hypoelliptic Laplacian is thus $\Delta_H=X_1^2+X_2^2$, with the $X_i$ given by \r{e:Xi-G5}.

We thus want to find the fundamental solution for the following heat equation:
\bqn
\partial_t \phi(t,x)=\Dh \phi(t,x),
\eqnl{e:hypo-G5}
with $\Delta_H=X_1^2+X_2^2$.
\subsubsection{Representations of $\GG_5$}

\renewcommand{\Rep}{\Repbase^{\lam,\mu,\nu}}
\renewcommand{\dRep}[1]{\dRepbase_{#1}^{\lam,\mu,\nu}}
\renewcommand{\domR}{\domRbase}
\renewcommand{\fzR}{\fzRbase}
\renewcommand{\DHL}{\DHLbase^{\lam,\mu,\nu}}

We now recall the representations of the group $\GG_5$, as computed by Dixmier in \cite[p. 338]{dixmier3}. As stated before, we may consider only representations  on the support of the Plancherel measure.
\bp
The dual space of $\GG_5$ is $\hat{G}=\Pg{\Rep\ |\ \lam^2+\mu^2\neq 0,\nu\in\R}$, where
\mmfunz{\Rep(x_1,x_2,x_3,x_4,x_5)}{\domR}{\domR}{\fzR(\th)}{\exp{i K^{\lam,\mu,\nu}_{x_1,x_2,x_3,x_4,x_5}(\th)}
\fzR\Pt{\th+\frac{\lam x_1+\mu x_2}{\lam^2+\mu^2}}}
with 
\bqn K^{\lam,\mu,\nu}_{x_1,x_2,x_3,x_4,x_5}(\th)&=&
-\frac12 \frac{\nu}{\lam^2+\mu^2}\Pt{\mu x_1 -\lam x_2}
+\lam x_4 +\mu x_5+\nn
&-&\frac16 \frac{\mu}{\lam^2+\mu^2}\Pt{\lam^2 x_1^3+3 \lam \mu x_1^2 x_2 + 3 \mu^2 x_1 x_2^2-\lam \mu x_2^3}
+\mu^2 x_1 x_2 \th
+\lam \mu (x_1^2-x_2^2)\th+\nn
&+&\frac12\Pt{\lam^2+\mu^2}\Pt{\mu x_1-\lam x_2} \th^2.\eqnn
The domain of $\Rep(x_1,x_2,x_3,x_4,x_5)$ is $\domR=L^2(\R,\C)$, endowed with the standard product $<\fzR_1,\fzR_2>:=\int_\R \fzR_1(\th)\overline{\fzR_2}(\th)\,d\th$ where $d\th$ is the Lebesgue measure.

The Plancherel measure on $\hat{G}$ is $dP(\lam,\mu,\nu)=d\lam d\mu d\nu$, i.e. the Lebesgue measure on $\R^3$.
\ep
\brem
Notice that in this case the domain $\domR$ of the representation $\Rep$ does not depend on $\lam,\mu,\nu$.
\erem

\subsubsection{The kernel of the hypoelliptic heat equation}

Consider the representation $\Rep$ of $\GG_5$ and let $\dRep{i}$ be the corresponding representations of the differential operators $L_{X_i}$ with $i=1,2$. Recall that $\dRep{i}$  are operators on $\domR$. Again from \cite[p. 338]{dixmier3}, or by explicit computation, we have
\bqn
\Pq{\dRep{1} \fzRbase}(\th)&=&
\Pt{-\frac{i}{2}\frac{\mu\eta}{\lam^2+\mu^2}
+\frac{\lam}{\lam^2+\mu^2}\frac{d}{d\th}
-\frac{i}2 \mu \Pt{\lam^2+\mu^2}\th^2}\fzRbase(\th)\nn
\Pq{\dRep{2}\fzRbase}(\th)&=&
\Pt{\frac{i}{2}\frac{\lam\mu}{\lam^2+\mu^2}
+\frac{\mu}{\lam^2+\mu^2}\frac{d}{d\th}
+\frac{i}2 \lam \Pt{\lam^2+\mu^2}\th^2}\fzRbase(\th),\eqnn
thus
\bqn
\Pq{\DHL\fzRbase}(\th)=\frac1{\lam^2+\mu^2}\frac{d^2 \fzRbase(\th)}{d\th^2}
-\frac{\Pt{\nu + (\lam^2 + \mu^2)^2 \th^2}^2}{4\Pt{\lam^2+\mu^2}}.\eqnn
The GFT of the hypoelliptic heat equation is thus
\bqn
\partial_t\fzRbase= \frac1{\lam^2+\mu^2}\frac{d^2 \fzRbase(\th)}{d\th^2}
-\frac{\Pt{\nu + (\lam^2 + \mu^2)^2 \th^2}^2}{4\Pt{\lam^2+\mu^2}}\fzRbase(\th).\eqnl{e:GFT-G5}
We rewrite it as
\bqn
\partial_\tau\fzRbase= \Pt{\frac{d^2}{d\th^2}-\Pt{\al\th^2+\beta}^2}\fzRbase(\th),\eqnl{e:GFT-G5-trasf}
with $\tau=\frac{t}{\Pt{\lam^2+\mu^2}}$, $\al=\frac{\lam^2+\mu^2}2$, $\beta=-\frac{\nu}{2}$.

The operator $\frac{d^2}{d\th^2}-(\al\th^2+\beta)^2$ is the Laplacian with quartic potential, see e.g. \cite{skala}. As already stated, no general explicit solutions are known for this equation. We call 
\bqn\sol_\tau \Pt{\th,\bar\th;\al,\beta}
\eqnn
the solution of
\bqn
\begin{cases}
\partial_\tau \psi(\tau,\th)=\Pt{\frac{d^2}{d\th^2}-\Pt{\al\th^2+\beta}^2}\psi(\tau,\th),\\
\psi(0,\th)=\de_{\bar\th},
\end{cases}
\eqnn
 i.e. the solution of \r{e:GFT-G5-trasf} evaluated in $\th$ at time $\tau$, with initial data $\de_{\bar\th}$ and parameters $\al,\beta$.

Applying Corollary \ref{c-2} and after straightforward computations, one gets  the kernel of the hypoelliptic heat equation on the group $\GG_5$:\\

$p_t(x_1,x_2,x_3,x_4,x_5)=$
\bqn =\int_{\lam^2+\mu^2\neq 0} d\lam d\mu d\nu \int_\R d \th\,
\exp{i K^{\lam,\mu,\nu}_{x_1,x_2,x_3,x_4,x_5}(\th)}
\sol_\frac{t}{\lam^2+\mu^2}\Pt{\th+\frac{\lam x_1+\mu x_2}{\lam^2+\mu^2},\th;\frac{\lam^2+\mu^2}2,-\frac{\nu}2}.
\eqnl{eq-G5-heat-exp}

%%%%%%%%%%%%%%%%%%%%%%%%%%%%% Biblio %%%%%%%%%%%%%%%%%%%%%%%%%%%%%%%%%%%%%5

\end{document}